\documentclass[12pt]{article}
\usepackage{amsmath}
\usepackage{amssymb}
\usepackage{amsthm}
\usepackage{amsfonts}
\usepackage{pstricks}
\usepackage{graphicx,psfrag}

      \newcommand {\al}   {\alpha}          \newcommand {\bt}  {\beta}
      \newcommand {\gam } {\gamma}

      \newcommand {\pl}   {\partial}        
           \newcommand {\RRR}  {{\mathbb R}}
           
     \newcommand {\EEE}  {{\mathcal E}}     \newcommand {\HHH}  {{\mathcal H}}
      
            \newcommand {\vk}   {\varkappa}

\newtheorem{theorem}{Theorem}[section]

\newcommand{\beq}{\begin{equation}}
\newcommand{\eeq}{\end{equation}}
\newtheorem{opr}{Definition}
\begin{document}

\title{Body with mirror surface and connected interior invisible from one point}

\author{ A. Aleksenko\thanks{Department of Mathematics, Aveiro University, Aveiro 3810, Portugal. This work was supported by Portuguese funds through the CIDMA - Center for Research and Development in Mathematics and Applications, and the Portuguese Foundation for Science and Technology ("FCT–Fund\c{c}\~{a}o para a Ci\^{e}ncia e a Tecnologia"), within project PEst-OE/MAT/UI4106/2014.
 }}

\maketitle

\begin{abstract}
 Here we demonstrate existence of a piecewise smooth obstacle having connected interior and invisible from a point in the framework of geometric optics.
\end{abstract}

The scattering theory prohibits existence of absolutely invisible bodies, since a nontrivial outgoing solution of the Helmholtz equation
cannot have zero scattering amplitude. Nevertheless, invisibility is possible in the framework of geometric optics which involves mathematical design of bodies with well-defined surfaces whose scattering map preserves certain trajectories of a flow of elastic particles. The main practical application of this study is optical shielding: by surrounding an object by a specially designed mirror surface, it is possible to create an illusion of invisibility from given points or directions.

The first work that targets the problem of designing a body invisible in a direction in the framework of mirror invisibility appears in \cite{0-resist} and is motivated by the problem of constructing a nonconvex body or zero resistance. The authors demonstrated that there exists a (connected and even simply connected) body invisible in one direction: if this body is manufactured out of perfectly reflective mirrors, a laser beam sent through this construction in the direction of invisibility would leave the body along the same trajectory. Remarkably, in \cite{Lak2012},\cite{Lak2012B},\cite{Lak2011} the scattering of acoustic waves by this body was studied.

This pioneering research led to several intriguing mathematical problems. One of them, proposed by Sergei Tabachnikov \cite{Tabach}, asks whether it is possible to design a body with mirror surface invisible in two directions. The problem was solved by Plakhov and Roshchina in \cite{PlakhovBook}: it was shown that a construction combining several pieces of parabolic cylinders can be used to produce a body invisible in two directions in the three-dimensional case. This body consists of two connected components,
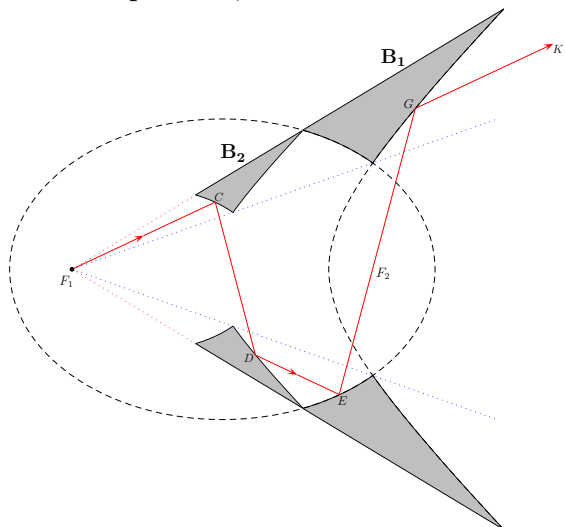
\begin{figure}[h]
\begin{picture}(0,170)

\rput(6.3,2.9){
\scalebox{2}{
\rput(-0.57,0){
\scalebox{0.534}{
\psdots[dotsize=0.7pt](-1,0)
\rput(-1.07,-0.15){\scalebox{0.4}{$F_1$}}
\rput(2.87,-0.05){\scalebox{0.4}{$F_2$}}
\rput(3,2.64){\scalebox{0.6}{$\mathbf{B_1}$}}
\rput(1,1.45){\scalebox{0.6}{$\mathbf{B_2}$}}
\pspolygon[linewidth=0.03pt,linecolor=white,fillstyle=solid,fillcolor=lightgray]
(0.534,0.925)(1.871,1.732)(1.64,1.483)(1.48,1.304)(1.34,1.14)(1.223,1)(1.18,0.947)
(1.136,0.89)(1.093,0.8367)(1.0456,0.7746)(1,0.709)(0.9,0.773)(0.8,0.826)(0.7,0.87)(0.6,0.906)
\pspolygon[linewidth=0.06pt,linecolor=white,fillstyle=solid,fillcolor=lightgray]
(0.534,-0.925)(1.871,-1.732)(1.64,-1.483)(1.48,-1.304)(1.34,-1.14)(1.223,-1)(1.18,-0.947)
(1.136,-0.89)(1.093,-0.8367)(1.0456,-0.7746)(1,-0.709)(0.9,-0.773)(0.8,-0.826)(0.7,-0.87)(0.6,-0.906)
\rput(0.83,0.9){\scalebox{0.35}{$C$}}
\rput(1.2,-1.1){\scalebox{0.35}{$D$}}
\rput(5.05,2.75){\scalebox{0.35}{$K$}}
\psline[linewidth=0.2pt,linestyle=dotted,linecolor=red,dotsep=1.6pt](-1,0)(0.534,-0.925)
\psline[linewidth=0.2pt,linestyle=dotted,linecolor=red,dotsep=1.6pt](-1,0)(0.534,0.925)
\psecurve[linewidth=0.12pt](0.45,0.948)(0.534,0.925)(0.6,0.9055)(0.7,0.8689)(0.8,0.8246)(0.9,0.7714)(1,0.7071)(1.1,0.6285)
\psecurve[linewidth=0.12pt](0.45,-0.948)(0.534,-0.925)(0.6,-0.9055)(0.7,-0.8689)
(0.8,-0.8246)(0.9,-0.7714)(1,-0.7071)(1.1,-0.6285)
\pscurve[linewidth=0.12pt]
(1.871,-1.732)(1.643,-1.483)(1.483,-1.304)(1.342,-1.14)
(1.2247,-1)(1.1832,-0.9487)(1.1402,-0.8944)(1.0954,-0.8367)(1.0488,-0.7746)(1,-0.7071)
\pscurve[linewidth=0.12pt]
(1,0.7071)(1.048,0.7746)(1.0954,0.8367)(1.1402,0.8944)(1.1832,0.9487)(1.2247,1)
(1.342,1.14)(1.483,1.304)(1.643,1.483)(1.871,1.732)
\psline[linewidth=0.12pt](0.534,0.925)(1.871,1.732)
\psline[linewidth=0.12pt](0.534,-0.925)(1.871,-1.732)
}}
\psdots[dotsize=0.15pt](-1,0)
\pspolygon[linewidth=0.03pt,linecolor=white,fillstyle=solid,fillcolor=lightgray]
(0.534,0.925)(1.871,1.732)(1.64,1.483)(1.48,1.304)(1.34,1.14)(1.223,1)(1.18,0.947)
(1.136,0.89)(1.093,0.8367)(1.0456,0.7746)(1,0.709)(0.9,0.773)(0.8,0.826)(0.7,0.87)(0.6,0.906)
\pspolygon[linewidth=0.03pt,linecolor=white,fillstyle=solid,fillcolor=lightgray]
(0.534,-0.925)(1.871,-1.732)(1.64,-1.483)(1.48,-1.304)(1.34,-1.14)(1.223,-1)(1.18,-0.947)
(1.136,-0.89)(1.093,-0.8367)(1.0456,-0.7746)(1,-0.709)(0.9,-0.773)(0.8,-0.826)(0.7,-0.87)(0.6,-0.906)
\psellipse[linewidth=0.01pt,linestyle=dashed,dash=1.2pt 0.8pt](0,0)(1.414,1)
\pscurve[linewidth=0.01pt,linestyle=dashed,dash=1.2pt 0.8pt]
(1.871,-1.732)(1.643,-1.483)(1.483,-1.304)(1.342,-1.14)
(1.2247,-1)(1.1832,-0.9487)(1.1402,-0.8944)(1.0954,-0.8367)(1.0488,-0.7746)(1,-0.7071)(0.9487,-0.6325)
(0.8944,-0.5477)(0.8367,-0.4472)(0.7746,-0.3162)(0.7071,0)(0.7746,0.3162)(0.8367,0.4472)(0.8944,0.5477)
(0.9487,0.6325)(1,0.7071)(1.048,0.7746)(1.0954,0.8367)(1.1402,0.8944)(1.1832,0.9487)(1.2247,1)
(1.342,1.14)(1.483,1.304)(1.643,1.483)(1.871,1.732)
\psline[linewidth=0.15pt,linestyle=dotted,dotsep=0.8pt,linecolor=blue](-1,0)(1.8284,1)
\psline[linewidth=0.15pt,linestyle=dotted,dotsep=0.8pt,linecolor=blue](-1,0)(1.8284,-1)
\pscurve[linewidth=0.06pt]
(1.871,-1.732)(1.643,-1.483)(1.483,-1.304)(1.342,-1.14)
(1.2247,-1)(1.1832,-0.9487)(1.1402,-0.8944)(1.0954,-0.8367)(1.0488,-0.7746)(1,-0.7071)
\pscurve[linewidth=0.06pt]
(1,0.7071)(1.048,0.7746)(1.0954,0.8367)(1.1402,0.8944)(1.1832,0.9487)(1.2247,1)
(1.342,1.14)(1.483,1.304)(1.643,1.483)(1.871,1.732)
\psecurve[linewidth=0.06pt](0.45,0.948)(0.534,0.925)(0.6,0.9055)(0.7,0.8689)(0.8,0.8246)(0.9,0.7714)(1,0.7071)(1.1,0.6285)
\psecurve[linewidth=0.06pt](0.45,-0.948)(0.534,-0.925)(0.6,-0.9055)(0.7,-0.8689)
(0.8,-0.8246)(0.9,-0.7714)(1,-0.7071)(1.1,-0.6285)
\psline[linewidth=0.06pt](0.534,0.925)(1.871,1.732)
\psline[linewidth=0.06pt](0.534,-0.925)(1.871,-1.732)
\psline[linewidth=0.15pt,linecolor=red,arrows=->,arrowscale=0.7](-1,0)(-0.5235,0.2235)
\psline[linewidth=0.15pt,linecolor=red,arrows=->,arrowscale=0.7]
(-0.5235,0.2235)(-0.047,0.447)(0.22,-0.571)(0.4985,-0.702)
\psline[linewidth=0.15pt,linecolor=red,arrows=->,arrowscale=0.7]
(0.4985,-0.702)(0.777,-0.833)(1.283,1.07)(2.2,1.5)
\psdots[dotsize=0.9pt](-1,0)
\rput(0.8,-0.87){\scalebox{0.2}{$E$}}
\rput(1.24,1.1){\scalebox{0.2}{$G$}}
}}
\end{picture}
\caption{A body invisible from one point.}
\label{fig invisible from 1 point}
\end{figure}
and its interior consists of 8 connected components, so it looks complicated to use such a construction in practical applications.
The main result of the paper is the following Theorem 1 (see figure \ref{fig invis 2D}).
\begin{theorem} Given a point in $\mathbb R^3$, there exists a body in $\mathbb R^3$ with connected interior which is
invisible from this point.
\end{theorem}
\section{Definitions}

We begin with reminding relevant definitions, then explain our construction and prove that it is invisible from two points.

\begin{opr}\label{def3}\rm
A {\it body} is a finite or countable union of its connected components, where each component is an open bounded domain with piecewise smooth boundary.
\end{opr}

\begin{opr}\label{def1}\rm
A body $B \subset \RRR^d$ is said to be {\it invisible from a point} $O \in \RRR^d \setminus B$, if for almost all $v \in S^{d-1}$ the billiard particle in $\RRR^d \setminus B$ emanating from $O$ with the initial velocity $v$, after a finite number of reflections from $\pl B$ will eventually move freely with the same velocity $v$ along a straight line containing $O$.
\end{opr}

If the point $O$ is infinitely distant, we get the notion of a body invisible in a direction. 


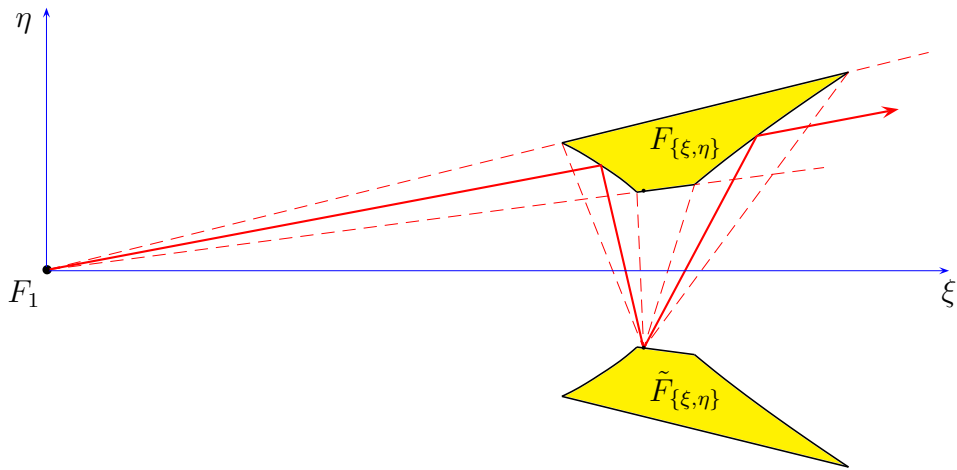
\begin{figure}[h]
\begin{picture}(0,183)
\rput{0}(-2.7,11){
\scalebox{1}{
\rput{-7.51}(3.5,-8){
\scalebox{1}{
\pspolygon[linecolor=white,fillstyle=solid,fillcolor=yellow]
(8.3863,2.25)(8.6207,2.5)(8.8667,2.75)(9.1225,3)(9.3868,3.25)(9.6584,3.5)(9.9363,3.75)(10.2197,4)
(6.565,2.57)(6.7015,2.5244)(6.9084,2.4402)(7.12,2.3415)(7.3,2.254)(7.4835,2.1521)(7.64,2.05)
\psdots[dotsize=3.5pt](0,0)
\pscurve[linewidth=0.6pt](7.64,2.05)(7.4835,2.1521)(7.3,2.254)(7.12,2.3415)(6.9084,2.4402)(6.7015,2.5244)(6.565,2.57)
\pscurve[linewidth=0.6pt](8.3863,2.25)(8.6207,2.5)(8.8667,2.75)(9.1225,3)(9.3868,3.25)(9.6584,3.5)(9.9363,3.75)(10.2197,4)
\psline[linewidth=0.2pt,linecolor=red,linestyle=dashed]
(0,0)(6.565,2.57)(8,0)(10.2197,4)
(11.2417,4.4)
\psline[linewidth=0.2pt,linecolor=red,linestyle=dashed]
(0,0)(7.64,2.05)(8,0)(8.3863,2.25)
(10.0636,2.7)
\psline[linewidth=0.6pt](6.565,2.57)(10.2197,4)
\psline[linewidth=0.6pt](7.64,2.05)(8.3863,2.25)
\psline[linewidth=0.8pt,linecolor=red,arrows=->,arrowscale=1.5](0,0)(7.12,2.3415)(8,0)(9.1225,3)(10.947,3.6)
\psdots[dotsize=1pt](8,0)
}}
\rput(3.56,-8.03){
\psline[linewidth=0.3pt,linecolor=blue,arrows=<->,arrowscale=1.3](12,0)(0,0)(0,3.5)
\rput(12,-0.3){$\xi$}
\rput(-0.3,3.3){$\eta$}
\rput(8.5,1.7){$F_{\{\xi,\eta\}}$}
\rput(-0.3,-0.3){$F_1$}
}
\rput{7.51}(3.5,-8.03){
\scalebox{1}{
\pspolygon[linecolor=white,fillstyle=solid,fillcolor=yellow]
(8.3863,-2.25)(8.6207,-2.5)(8.8667,-2.75)(9.1225,-3)(9.3868,-3.25)(9.6584,-3.5)(9.9363,-3.75)(10.2197,-4)
(6.565,-2.57)(6.7015,-2.5244)(6.9084,-2.4402)(7.12,-2.3415)(7.3,-2.254)(7.4835,-2.1521)(7.64,-2.05)
\psdots[dotsize=1pt](0,0)(8,0)
\pscurve[linewidth=0.6pt](7.64,-2.05)(7.4835,-2.1521)(7.3,-2.254)(7.12,-2.3415)(6.9084,-2.4402)(6.7015,-2.5244)
(6.565,-2.57)
\pscurve[linewidth=0.6pt](8.3863,-2.25)(8.6207,-2.5)(8.8667,-2.75)(9.1225,-3)(9.3868,-3.25)(9.6584,-3.5)
(9.9363,-3.75)(10.2197,-4)
\psdots[dotsize=2pt](0,0)
\psdots[dotsize=1.5pt](8,0)(7.7274,-2.0706)
\psline[linewidth=0.6pt](6.565,-2.57)(10.2197,-4)
\psline[linewidth=0.6pt](7.64,-2.05)(8.3863,-2.25)
}}
\rput(3.56,-8.03){
\rput(8.5,-1.6){$\tilde F_{\{\xi,\eta\}}$}
}
}}
\end{picture}
\caption{A two-dimensional figure invisible from the origin. The 3-dimensional construction is obtained by rotating this figure around the $\xi$-axis.}
\label{fig invis 2D}
\end{figure}

Notice that a 3D body invisible from one point was constructed in \cite{PlakhovBook} (a central cross section of this body by a plane passing through the point is shown in Fig.~\ref{fig invis 2D}). Its interior is disconnected: it consists two connected components. This provides a difficulty in practical realization of this construction. On the contrary, below we construct a body with connected interior.

\section{Construction}
We describe the geometrical shape of the body invisible from one point, provide a proof of its invisibility, and then give exact formulas that determine its shape. The description is made in several steps.
\vspace{2mm}

{\bf 1.} Consider the ellipse $\EEE$ given by
$$
\frac{x^2}{a^2} + \frac{y^2}{b^2} = 1, \quad a > 0, \, \ b > 0
$$
in Cartesian coordinates $x,\, y$. The foci of $\EEE$ are the points $F_1 = (-c, 0)$ and $F_2 = (c, 0)$, where $c = \sqrt{a^2 - b^2}$. Next consider the hyperbola
$$
\frac{x^2}{\al^2} - \frac{y^2}{\bt^2} = 1, \quad \al > 0, \, \ \bt > 0.
$$
We require that the hyperbola has the same foci $F_1$ and $F_2$, that is, the parameters $\al$ and $\bt$ satisfy the equality
\beq\label{c al bt}
c = \sqrt{\al^2 + \bt^2}.
\eeq

\begin{figure}[h]
\begin{picture}(0,210)
\rput(6.0,3.0){
\scalebox{2.3}{
\psdots[dotsize=0.15pt](-1,0)(0.534,0.925)(1.871,1.732)(1.241,0.479)
\pspolygon[linewidth=0.03pt,linecolor=white]
(0.534,0.925)(1.871,1.732)(1.64,1.483)(1.48,1.304)(1.34,1.14)(1.223,1)(1.18,0.947)
(1.136,0.89)(1.093,0.8367)(1.0456,0.7746)(1,0.709)(0.9,0.773)(0.8,0.826)(0.7,0.87)(0.6,0.906)
\pspolygon[linewidth=0.03pt,linecolor=white]
(0.534,-0.925)(1.871,-1.732)(1.64,-1.483)(1.48,-1.304)(1.34,-1.14)(1.223,-1)(1.18,-0.947)
(1.136,-0.89)(1.093,-0.8367)(1.0456,-0.7746)(1,-0.709)(0.9,-0.773)(0.8,-0.826)(0.7,-0.87)(0.6,-0.906)
\psellipse[linewidth=0.233pt](0,0)(1.414,1)
\psdots[dotsize=0.7pt](-1,0)(1,0)
\rput(-1.05,-0.07){\scalebox{0.25}{$F_1$}}
\rput(1.07,-0.01){\scalebox{0.25}{$F_2$}}
  \psarc[linewidth=0.2pt,linecolor=brown](1,0){0.145}{63}{90}
  \psarc[linewidth=0.2pt,linecolor=brown](1,0){0.11}{90}{117}
  \psarc[linewidth=0.2pt,linecolor=brown](1,0){0.13}{63}{90}
  \psarc[linewidth=0.2pt,linecolor=brown](1,0){0.095}{90}{117}
  \rput(0.95,0.2){\scalebox{0.3}{$\al$}}
  \rput(1.06,0.23){\scalebox{0.3}{$\bt$}}
  \rput(-0.45,1.05){\scalebox{0.45}{$\EEE$}}
  \rput(1.5,-1.1){\scalebox{0.45}{$\HHH$}}
\pscurve[linewidth=0.233pt]
(1.483,-1.304)(1.342,-1.14)
(1.2247,-1)(1.1832,-0.9487)(1.1402,-0.8944)(1.0954,-0.8367)(1.0488,-0.7746)(1,-0.7071)(0.9487,-0.6325)
(0.8944,-0.5477)(0.8367,-0.4472)(0.7746,-0.3162)(0.7071,0)(0.7746,0.3162)(0.8367,0.4472)(0.8944,0.5477)
(0.9487,0.6325)(1,0.7071)(1.048,0.7746)(1.0954,0.8367)(1.1402,0.8944)(1.1832,0.9487)(1.2247,1)
(1.342,1.14)(1.483,1.304)(1.643,1.483)(1.871,1.732)
\rput(0.91,0.68){\scalebox{0.25}{$C$}}
\rput(0.52,1){\scalebox{0.25}{$A$}}
\rput(1.85,1.78){\scalebox{0.25}{$B$}}
  \psline[linewidth=0.15pt,linestyle=dashed,dash=1.2pt 0.8pt,linecolor=blue](0.534,0.925)(1,0)
  (1.871,1.732)(-1,0)
\psline[linewidth=0.15pt,linecolor=red](-1,0)(1,0)(1,0.7071)
}}
\end{picture}
\caption{Ellipse and hyperbola.}
\label{fig ellipse and parabola}
\end{figure}
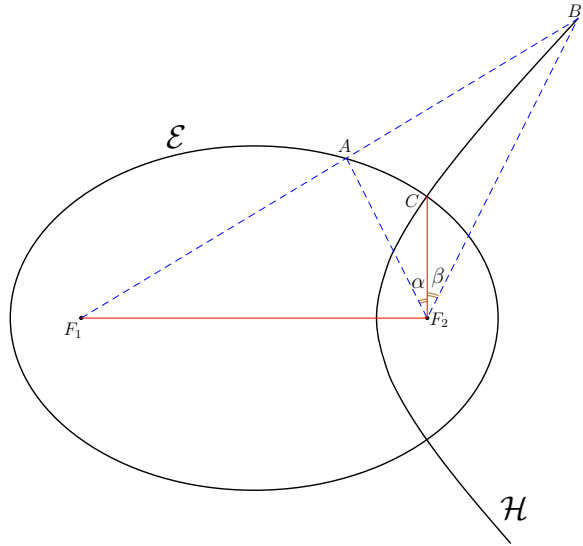

Denote by $\HHH$ the right branch of the hyperbola. There are two points of intersection of $\HHH$ with the ellipse $\EEE$, which are symmetrical to each other with respect to the $x$-axis; we denote by $C$ the upper point of intersection (see Fig. \ref{fig ellipse and parabola}).
Let us additionally impose the condition that the segments $F_1 F_2$ and $F_2 C$ are perpendicular; one easily sees that this condition is equivalent to the equation
$$
\al a = c^2.
$$
It is convenient to introduce the parameter
\beq\label{kappa}
\vk = \frac{a}{c} = \frac{c}{\al}.
\eeq
\vspace{2mm}

{\bf 2.} Here we prove some auxiliary geometric statements which will be needed later on. First state a characteristic property of angle bisector in a triangle.
\vspace{1mm}

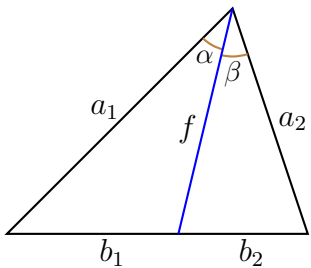
\begin{figure}[h]
\begin{picture}(0,110)
\rput(4.5,1){
\scalebox{1}{
\pspolygon(0,0)(4,0)(3,3)
\psarc[linecolor=brown](3,3){0.56}{225.3}{256.5}
\psarc[linecolor=brown](3,3){0.64}{256.7}{288.2}
\rput(2.63,2.36){\scalebox{0.9}{$\al$}}
\rput(3,2.13){\scalebox{0.9}{$\bt$}}
\psline[linewidth=0.8pt,linecolor=blue](3,3)(2.28,0)
\rput(1.3,1.65){$a_1$}
\rput(3.8,1.5){$a_2$}
\rput(1.4,-0.25){$b_1$}
\rput(3.25,-0.25){$b_2$}
\rput(2.4,1.4){$f$}
}}
\end{picture}
\caption{The characteristic property of the angle bisector.
}
\label{fig bisector}
\end{figure}

{\bf Property.} {\it The segment $f$ is the bisector of the corresponding angle in Figure \ref{fig bisector} (that is, $\al = \bt$), if and only if $(a_1 + b_1)(a_2 - b_2) = f^2$.}
\vspace{1mm}

{\it Sketch of the proof.} Consider the following relations on the values $a_1$,\, $a_2$,\, $b_1$,\, $b_2$, and $f$:
$$
1. \quad a_1/a_2 = b_1/b_2; \hspace*{91mm}
$$
$$
2. \quad a_1 a_2 - b_1 b_2 = f^2; \hspace*{85mm}
$$
\beq\label{bisector property}
3. \quad (a_1 + b_1)(a_2 - b_2) = f^2. \hspace*{67mm}
\eeq
The equalities 1 and 2 are well known in the literature; each of them is a characteristic property of triangle bisector. The equality 3 is a direct consequence of the equalities 1 and 2; thus the direct property (\ref{bisector property}) of the angle bisector is established. The proof of the inverse property (\ref{bisector property}) is also simple, but cumbersome, and utilizes the sine rule and some trigonometry. It is omitted here. \hspace{95mm} $\Box$
\vspace{1mm}

{\bf Proposition.} {\it The angles $\al = \measuredangle A F_2 C$ and $\bt = \measuredangle B F_2 C$ in Figure \ref{fig ellipse and parabola} are equal.}
\vspace{1mm}

\begin{proof}
Let us make an auxiliary construction. Extend the segment $BF_2$ until the second intersection with the ellipse at a point $A'$. Denote by $C'$ the second point of intersection of the ellipse with the branch of the hyperbola $\HHH$. Denote
$$
f = 2c = |F_1 F_2|, \ \, g = |F_2 C| = |F_2 C'|, \ \, a_1 = |F_1 A'|,
$$
$$
b_1 = |F_2 A'|, \ \, a_2 = |F_1 B|, \ \text{and} \ b_2 = |F_2 B|
$$
(see Fig. \ref{fig proof}).
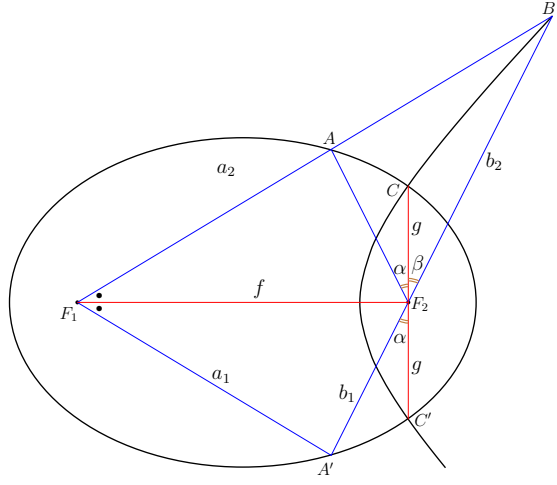
\begin{figure}[h]
\begin{picture}(0,200)
\rput(6.0,3.2){
\scalebox{2.2}{
\psdots[dotsize=0.15pt](-1,0)(0.534,0.925)(1.871,1.732)(1.241,0.479)
\pspolygon[linewidth=0.03pt,linecolor=white]
(0.534,0.925)(1.871,1.732)(1.64,1.483)(1.48,1.304)(1.34,1.14)(1.223,1)(1.18,0.947)
(1.136,0.89)(1.093,0.8367)(1.0456,0.7746)(1,0.709)(0.9,0.773)(0.8,0.826)(0.7,0.87)(0.6,0.906)
\pspolygon[linewidth=0.03pt,linecolor=white]
(0.534,-0.925)(1.871,-1.732)(1.64,-1.483)(1.48,-1.304)(1.34,-1.14)(1.223,-1)(1.18,-0.947)
(1.136,-0.89)(1.093,-0.8367)(1.0456,-0.7746)(1,-0.709)(0.9,-0.773)(0.8,-0.826)(0.7,-0.87)(0.6,-0.906)
\psellipse[linewidth=0.233pt](0,0)(1.414,1)
\psdots[dotsize=0.7pt](-1,0)(1,0)
\rput(-1.05,-0.07){\scalebox{0.25}{$F_1$}}
\rput(1.07,-0.01){\scalebox{0.25}{$F_2$}}
\rput(-0.87,0.04){\scalebox{0.2}{$\bullet$}}
\rput(-0.87,-0.04){\scalebox{0.2}{$\bullet$}}
\pscurve[linewidth=0.233pt]
(1.2247,-1)(1.1832,-0.9487)(1.1402,-0.8944)(1.0954,-0.8367)(1.0488,-0.7746)(1,-0.7071)(0.9487,-0.6325)
(0.8944,-0.5477)(0.8367,-0.4472)(0.7746,-0.3162)(0.7071,0)(0.7746,0.3162)(0.8367,0.4472)(0.8944,0.5477)
(0.9487,0.6325)(1,0.7071)(1.048,0.7746)(1.0954,0.8367)(1.1402,0.8944)(1.1832,0.9487)(1.2247,1)
(1.342,1.14)(1.483,1.304)(1.643,1.483)(1.871,1.732)
\rput(0.91,0.68){\scalebox{0.25}{$C$}}
  \rput(1.09,-0.706){\scalebox{0.25}{$C'$}}
\rput(0.52,1){\scalebox{0.25}{$A$}}
  \rput(0.5,-1.01){\scalebox{0.25}{$A'$}}
\rput(1.85,1.78){\scalebox{0.25}{$B$}}
  \psline[linewidth=0.15pt,linecolor=blue](0.534,0.925)(1,0)
  \psline[linewidth=0.15pt,linecolor=blue](1.871,1.732)(-1,0)
  \psline[linewidth=0.15pt,linecolor=blue](-1,0)(0.534,-0.925)(1.871,1.732)
\rput(0.1,0.08){\scalebox{0.3}{$f$}}
\rput(1.05,0.45){\scalebox{0.3}{$g$}}
\rput(1.05,-0.39){\scalebox{0.3}{$g$}}
\rput(-0.13,-0.45){\scalebox{0.3}{$a_1$}}
\rput(0.63,-0.55){\scalebox{0.3}{$b_1$}}
\rput(1.52,0.85){\scalebox{0.3}{$b_2$}}
\rput(-0.1,0.8){\scalebox{0.3}{$a_2$}}
\psline[linewidth=0.15pt,linecolor=red](-1,0)(1,0)
\psline[linewidth=0.15pt,linecolor=red](1,-0.7071)(1,0.7071)
  \psarc[linewidth=0.2pt,linecolor=brown](1,0){0.145}{63}{90}
  \psarc[linewidth=0.2pt,linecolor=brown](1,0){0.11}{90}{117}
  \psarc[linewidth=0.2pt,linecolor=brown](1,0){0.13}{63}{90}
  \psarc[linewidth=0.2pt,linecolor=brown](1,0){0.095}{90}{117}
  \rput(0.95,0.2){\scalebox{0.3}{$\al$}}
  \rput(0.95,-0.22){\scalebox{0.3}{$\al$}}
    \psarc[linewidth=0.2pt,linecolor=brown](1,0){0.125}{243}{270}
  \psarc[linewidth=0.2pt,linecolor=brown](1,0){0.11}{243}{270}
  \rput(1.06,0.23){\scalebox{0.3}{$\bt$}}
}}
\end{picture}
\caption{Auxiliary construction.}
\label{fig proof}
\end{figure}
By the focal property of the ellipse, we have $|F_1 A'| + |F_2 A'| = |F_1 C'| + |F_2 C'|$, that is,
\beq\label{eq1}
a_1 + b_1 = \sqrt{f^2 + g^2} + g.
\eeq
Further, by the focal property of the hyperbola we have $|F_1 B| - |F_2 B| = |F_1 C| - |F_2 C|$, that is,
\beq\label{eq2}
a_2 - b_2 = \sqrt{f^2 + g^2} - g.
\eeq
Multiplying both sides of (\ref{eq1}) and (\ref{eq2}), we get
$$
(a_1 + b_1)(a_2 - b_2) = f^2,
$$
and taking into account the Property, one concludes that $F_1 F_2$ is the bisector of the angle $F_1$ in the triangle $A'F_1B$. This means that  $A'$ is symmetric to $A$ with respect to the straight line $F_1 F_2$, and by symmetry one has
\beq\label{eq3}
\measuredangle AF_2C = \measuredangle A'F_2C'.
\eeq
On the other hand, the angles $\measuredangle BF_2C$ and $\measuredangle A'F_2C'$ are vertical, and therefore, are equal:
\beq\label{eq4}
\measuredangle BF_2C = \measuredangle A'F_2C'.
\eeq
The equations (\ref{eq3}) and (\ref{eq4}) imply that $\measuredangle AF_2C = \measuredangle BF_2C$, therefore $\al = \bt$.
\end{proof}
\vspace{2mm}

{\bf 3.} Draw a ray with the vertex at $F_1$,
$$
y = k(x + c), \quad x \ge -c,
$$
with $k > 0$. The ray intersects the branch $\HHH$ of the hyperbola, if and only if $k < {\bt}/{\al}$. Taking into account the relations (\ref{c al bt}) and (\ref{kappa}) on $\al$ and $\bt$, one rewrites this inequality as $k < k_{\max}$, where
\beq\label{k max}
k_{\max} = \sqrt{\vk^2 - 1}.
\eeq
Suppose that $k$ satisfies (\ref{k max}) and denote by $A$ and $B$ the points of intersection of the ray with $\EEE$ and $\HHH$, respectively (see Fig. \ref{fig ellipse and parabola}).

In what follows we will also assume that the inequalities
\beq\label{ineq dlin}
|F_1 A| < |F_1 F_2| < |F_1 B|
\eeq
are satisfied. Below we derive the condition on $k$ equivalent to (\ref{ineq dlin}). Denote $A = (x_A, y_A)$ and $B = (x_B, y_B)$; the following relations can be easily derived:
\beq\label{a relation}
|F_1 A| = \frac ca\, x_A + a \quad \text{and} \quad |F_1 B| = \frac c\al\, x_B + \al.
\eeq
By the second formula in (\ref{a relation}), one has $|F_1 B| > |F_1 C| > |F_1 F_2|$, and so, the second inequality in (\ref{ineq dlin}) is always satisfied.

Note that
\beq\label{2c}
|F_1 F_2| = 2c.
\eeq
The ray with the largest inclination $y = k_{\max}(x + c)$ intersects $\EEE$ at the point $A_{\infty} = (0, b)$, therefore $|F_1 A_{\infty}| = \sqrt{c^2 + b^2} = a$. We impose the condition
$$
\vk < 2;
$$
then the distance $|F_1 A|$ monotonically decreases from $|F_1 C| = \sqrt{(2c)^2 + b^4/a^2} > 2c$ to $|F_1 A_{\infty}| = a < 2c$ when $A$ runs the elliptic curve $CA_{\infty}$ from $C$ to $A_{\infty}$, and takes the value $2c$ at a single point $A_0$ in between.

Using (\ref{2c}) and the first formula in (\ref{a relation}), we conclude that the first inequality in (\ref{ineq dlin}) is equivalent to $(c/a) x_A + a < 2c$, which can be rewritten as
$$
x_A < x_0 = a \left(2 - \frac ac\right).
$$
Let $A_0 = (x_0, y_0)$ be the point on the ellipse; then one has
$$
y_0 = c \sqrt{\vk^2 - 1} \sqrt{1 - (2 - \vk)^2}.
$$
We conclude that  the first inequality in (\ref{ineq dlin}) is equivalent to $k > k_{\min}$, where
\beq\label{k min}
k_{\min} = \frac{y_0}{x_0 + c} = \frac{\sqrt{\vk^2 - 1}\sqrt{1 - (2 - \vk)^2}}{1 + 2\vk - \vk^2} = (\vk - 1) \frac{\sqrt{4 - (\vk - 1)^2}}{2 - (\vk - 1)^2}.
\eeq

Thus, the condition ensuring that the ray $y = k(x + c), \ x \ge -c$ intersects both $\EEE$ and $\HHH$ and that for the points of intersection, $A$ and $B$, the inequalities (\ref{ineq dlin}) are satisfied, reads as
$$
k_{\min} < k < k_{\max}.
$$
\vspace{2mm}

\begin{figure}[h]
\begin{picture}(0,210)
\rput(1.5,3.5){
\scalebox{0.8}{
\pspolygon[linecolor=white,fillstyle=solid,fillcolor=yellow]
(8.3863,2.25)(8.6207,2.5)(8.8667,2.75)(9.1225,3)(9.3868,3.25)(9.6584,3.5)(9.9363,3.75)(10.2197,4)
(6.565,2.57)(6.7015,2.5244)(6.9084,2.4402)(7.12,2.3415)(7.3,2.254)(7.4835,2.1521)(7.64,2.05)
\psellipse[linewidth=0.8pt](4,0)(5,3)
\psdots[dotsize=1.3pt](0,0)(8,0)(8,1.8)(6.565,2.57)(8.3863,2.25)
\psdots[dotsize=3pt](0,0)(8,0)
\pscurve[linewidth=0.8pt](7.2,0)(7.2173,0.25)(7.2687,0.5)(7.3526,0.75)(7.4667,1)(7.608,1.25)(7.7736,1.5)(7.9604,1.75)
(8.1655,2)(8.3863,2.25)(8.6207,2.5)(8.8667,2.75)(9.1225,3)(9.3868,3.25)(9.6584,3.5)(9.9363,3.75)(10.2197,4)
\psline[linewidth=0.2pt,linecolor=red,linestyle=dashed]
(0,0)
(10.2197,4)(11.2417,4.4)
\psline[linewidth=0.2pt,linecolor=red,linestyle=dashed]
(0,0)
(8.3863,2.25)(10.0636,2.7)
\psline[linewidth=0.8pt](6.565,2.57)(10.2197,4)
\psline[linewidth=0.8pt](7.64,2.05)(8.3863,2.25)
\psline[linewidth=0.8pt,linecolor=red,arrows=->,arrowscale=1.5](0,0)(7.12,2.3415)(8,0)(9.1225,3)(10.947,3.6)
\psline[linewidth=0.8pt,linecolor=red,linestyle=dotted](7.12,2.3415)(9.1225,3)
\rput(7,2.1){\scalebox{0.8}{$\tilde A$}}
\rput(9.24,2.83){\scalebox{0.8}{$\tilde B$}}
\rput(7.5,1.88){\scalebox{0.7}{$A$}}
\rput(8.55,2.17){\scalebox{0.7}{$B$}}
\psline[linewidth=0.8pt](7.7,0)(8.3,0)
\rput(8.1,2.7){\scalebox{1}{$F_{\{x,y\}}$}}
\rput(-0.25,-0.25){\scalebox{0.8}{$F_1$}}
\rput(8.07,-0.3){\scalebox{0.8}{$F_2$}}
\rput(8.27,1.8){\scalebox{0.8}{$C$}}
\psline[linewidth=0.2pt,linestyle=dashed,linecolor=blue](0.2,0)(7.6,0)
\psline[linewidth=0.2pt,linestyle=dashed,linecolor=blue](8,0)(8,1.8)
\rput(11,3.34){\scalebox{0.8}{$D$}}
\rput(10,2.5){\scalebox{0.8}{$K$}}
}}
\end{picture}
\caption{A light ray reflecting from the mirrors.}
\label{fig 3 reflections}
\end{figure}
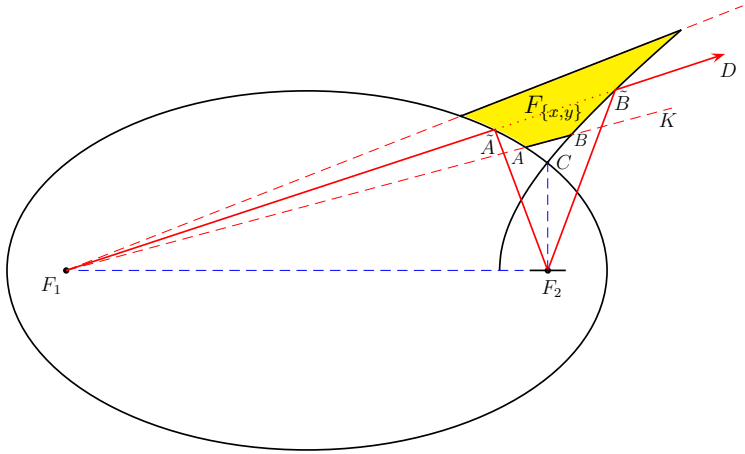

{\bf 4.} Draw two rays with inclinations $k_1$ and $k_2$,\, $y = k_1 (x + c),\, x \ge -c$ and $y = k_2 (x + c),\, x \ge -c$, where
\beq\label{k1k2}
k_{\min} < k_1 < k_2 < k_{\max}.
\eeq
The ray $y = k_1 (x + c),\, x \ge -c$ is denoted by $F_1 K$ in Figure \ref{fig 3 reflections}. From the previous item we know that both rays intersect $\EEE$ and $\HHH$ and the inequalities (\ref{ineq dlin}) are satisfied, with $A$ and $B$ being the points of intersection of $F_1 K$ with $\EEE$ and $\HHH$.

Determine the figure $F_{\{x,y\}}$ by
$$
\frac{x^2}{a^2} + \frac{y^2}{b^2} > 1, \quad \frac{x^2}{\al^2} - \frac{y^2}{\bt^2} < 1,
$$
$$
k_1 < \frac{y}{x + c} < k_2, \quad y > 0
$$
(see Fig. \ref{fig 3 reflections}).

Take a ray $F_1 D$ at an inclination $k \in (k_1,\, k_2)$. Let $\tilde A$ and $\tilde B$ be the points of intersection of this ray with the elliptic and hyperbolic arcs forming the boundary of $F_{\{x,y\}}$. Now imagine that the boundary of $F_{\{x,y\}}$ is mirror-like and there is a flat mirror on the line $F_1 F_2$. Then the broken line $F_1 \tilde A F_2 \tilde B D$  represents a light ray emanating from $F_1$ and making reflections from these mirror boundaries.

Indeed, according to the focal property of the billiard in ellipse, the light ray from $F_1$, after a reflection at $\tilde A$, gets into $F_2$. The segment $F_2 C$ is orthogonal to $F_1F_2$ and is the bisector of the angle $\tilde A F_2 \tilde B$, as proved in the Proposition. Therefore the light ray, after the second reflection at $F_2$, gets into $\tilde B$. According to the focal property of the billiard in hyperbola, the light ray reflected at $\tilde B$ moves along the straight line $\tilde B D$ through $F_1$.

Now take the angle $\gam = \frac 12 \arctan k_1 = \frac 12 \measuredangle KF_1F_2$. The tangent $t = \tan\gam$ satisfies the equation
\beq\label{t}
\frac{2t}{1 - t^2} = k_1,
\eeq
which implies that
$$
t = \frac{\sqrt{k_1^2 + 1} - 1}{k_1}.
$$
Make the change of variables
$$
\xi = \frac{(x + c) + ty}{\sqrt{1 + t^2}} = \cos\gam \cdot (x + c) + \sin\gam \cdot y,
$$
$$
\eta = \frac{-t(x + c) + y}{\sqrt{1 + t^2}} = -\sin\gam \cdot (x + c) + \cos\gam \cdot y.
$$
The inverse change of variables has the form
$$
x + c = \frac{\xi - t\eta}{\sqrt{1 + t^2}},
$$
$$
y = \frac{t\xi + \eta}{\sqrt{1 + t^2}}.
$$
The new coordinate system $\xi,\, \eta$ is orthogonal, its origin $\xi = 0,\, \eta = 0$ coincides with the point $F_1 = (-c, 0)$ (in the $x,y$-coordinates), and the $\xi$-axis (given by the equality $\eta = 0$) is the bisector of the angle $KF_1F_2$ formed by the lines $y = 0$ and $y = k_1(x + c)$.

In the new coordinates $\xi,\, \eta$ the figure $F_{\{x,y\}}$ takes the following form:
$$
F_{\{\xi,\eta\}} = \left\{ (\xi, \eta): \right. \hspace*{67mm}
$$
$$
\frac{(\xi - t\eta)^2}{\al^2} - \frac{(t\xi + \eta)^2}{\bt^2} < 1 + t^2 < \frac{(\xi - t\eta)^2}{a^2} + \frac{(t\xi + \eta)^2}{b^2},
$$
\beq\label{xi eta repr}
\hspace*{43mm} \left. k_1 < \frac{t\xi + \eta}{\xi - t\eta} < k_2, \quad t\xi + \eta > 0 \right\}.
\eeq
Let $\tilde F_{\{\xi,\eta\}}$ be symmetric to $F_{\{\xi,\eta\}}$ with respect to the line $\eta = 0$; then the two-dimensional figure $F_{\{\xi,\eta\}} \cup \tilde F_{\{\xi,\eta\}}$ is invisible from the origin $F_1$ (see Fig. \ref{fig invis 2D}).

Indeed, a light ray emanated from $F_1$ makes the first reflection from the elliptic arc bounding $F_{\{\xi,\eta\}}$. The second reflection is from a point on the flat segment bounding $\tilde F_{\{\xi,\eta\}}$, besides the distance from $F_1$ to this point equals $|F_1F_2|$. The condition (\ref{k1k2}) and the inequalities (\ref{ineq dlin}) ensure that this point really belongs to the flat segment.

\begin{figure}[h]
\centering
\includegraphics[scale=0.45]{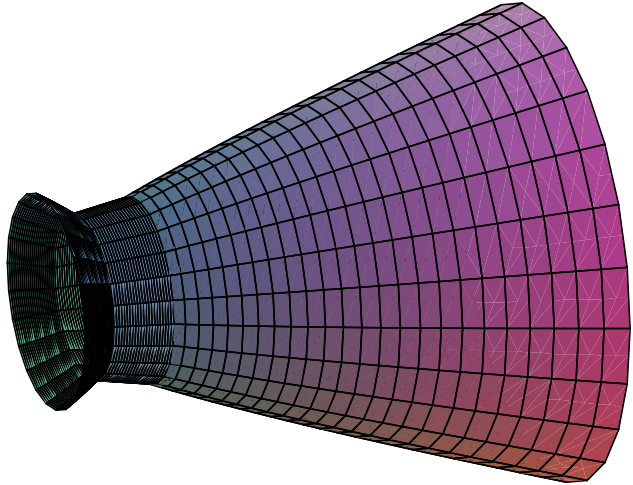}
\caption{The 3-dimensional body obtained by rotating the plane figure on Fig.~\ref{fig invis 2D} around the horizontal axis. In order to make the body's shape more visible, the exterior part of its boundary is removed.}
\label{fig:3D}
\end{figure}

The three-dimensional figures $G_1$ and $G_2$ invisible from the origin are obtained by rotating the figure $F_{\{\xi,\eta\}} \cup \tilde F_{\{\xi,\eta\}}$ with respect to the axis $\eta = 0$ and to the axis $\xi = 0$. In the first case (see Fig.~\ref{fig:3D}) the figure $G_1$ is
\beq\label{G1}
G_1 = \{ (u,v,w): (u, \sqrt{v^2+w^2}) \in F_{\{\xi,\eta\}} \};
\eeq
in the second case the figure $G_2$ is
\beq\label{G2}
G_2 = \{ (u,v,w): (\sqrt{u^2+v^2}, |w|) \in F_{\{\xi,\eta\}} \}.
\eeq

{\bf 5.}
Summarizing, the construction of an invisible body is as follows. Choose the parameters $c > 0$ and $1 < \vk < 2$. Calculate $k_{\min}$ and $k_{\max}$ according to the formulas (\ref{k min}) and (\ref{k max}), and choose the parameters $k_1$ and $k_2$ satisfying (\ref{k1k2}). Define $a^2,\, b^2,\, \al^2,\, \bt^2$ by
$$
a^2 = \vk^2 c^2, \quad b^2 = (\vk^2 - 1) c^2, \quad \al^2 = \vk^{-2} c^2, \quad \bt^2 = (1 - \vk^{-2}) c^2,
$$
and calculate $t$ according to (\ref{t}). Finally, define the 2D region $F_{\{\xi,\eta\}}$ by (\ref{xi eta repr}), and define the regions $G_1$ and $G_2$ in the three-dimensional space of Cartesian coordinates $u,v,w$ by (\ref{G1}) and (\ref{G2}). Each of these regions depends on 4 continuous parameters: scale of the picture $c$, excentricity of the ellipse $\vk$, and inclinations of two generating lines, $k_1$ and $k_2$.

\section*{Acknowledgements}
The author acknowledges Alexander Plakhov for the permanent support and advice.

\end{document}